\documentclass[10pt]{article}
\usepackage[matrix,arrow,curve]{xy}
\usepackage{amssymb,amsthm}
\usepackage{supertabular}
\usepackage[dvips]{graphics}

\newtheorem{dummy}{anything}
\newtheorem{Theorem}[dummy]{Theorem}
\newtheorem{lemma}[dummy]{Lemma}
\newtheorem{Proposition}[dummy]{Proposition}
\newtheorem{Corollary}[dummy]{Corollary}

\newtheorem{Example}[dummy]{Example}

\newcommand{\pcirc}{\kern .7pt {\scriptstyle \circ} \kern 1pt}

\newcommand{\R}{{\mathbb R}}
\newcommand{\C}{{\mathbb C}}
\newcommand{\Z}{{\mathbb Z}}

\newcommand{\vol}{{\rm {vol}}}

\newcommand{\eqref}[1]{(\ref{#1})}
\newcommand{\hfl}[2]{\smash{\mathop{\hbox to 1 truecm{\kern %
3pt\rightarrowfill\kern 3pt}}%
\limits^{\scriptstyle#1}_{\scriptstyle#2}}}
\newcommand{\cqfd}{\unskip\kern 6pt\penalty 500%
\raise -2pt\hbox{\vrule\vbox to10pt{\hrule width %
4pt\vfill\hrule}\vrule}\smallskip}

\title{Betti numbers of random manifolds}
\author{Michael Farber\\Department of Mathematical Sciences\\ University of Durham, UK
\and Thomas Kappeler\\ Institute of Mathematics\\ University of
Zurich, Switzerland}
\date{November 19, 2006}

\begin{document}
\maketitle 

\begin{abstract}
We study mathematical expectations of Betti numbers of
configuration spaces of planar linkages, viewing the lengths of
the bars of the linkage as random variables. Our main result gives
an explicit asymptotic formulae for these mathematical
expectations for two distinct probability measures describing the
statistics of the length vectors when the number of links tends to
infinity. In the proof we use a combination of geometric and
analytic tools. The average Betti numbers are expressed in terms
of volumes of intersections of a simplex with certain half-spaces.
\end{abstract}

\section{Introduction}

In various fields of applications, such as topological robotics,
configuration spaces of mechanical systems depend on a large
number of parameters, which typically are only partially known and
often can be considered as random variables. Since these
parameters determine the topology of the configuration space, the
latter can be viewed in such a case as a {\it random topological
space} or a {\it random manifold}. To control such a system one
has to understand geometry, topology and control theory of random
manifolds.

One of the most natural notion to investigate is the mathematical
expectation of the Betti numbers of random manifolds. Clearly,
these average Betti numbers encode valuable information for
engineering applications; for instance they provide an average
lower bound for the number of critical points of a Morse function
(i.e. observable) on such manifolds.

In this paper we consider a specific instance of this general
problem. We study closed planar $n$-gons whose sides have fixed
lengths $l_1, \dots, l_n$ where $l_{i}>0$ for $1 \le i \le n$. The
{\it \lq\lq polygon\rq\rq}\, space
\begin{eqnarray}
M_\ell = \{(u_1, \dots, u_n)\, \in \, S^1\times \dots\times S^1;
\, \sum_{i=1}^n l_iu_i \, =\, 0\in \C\}\, /\, {\rm {SO}}(2)
\end{eqnarray}
parametrizes the variety of all possible shapes of such planar
$n$-gons with sides of length $l_1, \dots, l_n$. The unit vector
$u_i\in \C$ indicates the direction of the $i$-th side of the
polygon. The condition $\sum l_iu_i=0$ expresses the property of
the polygon being closed. The rotation group ${\rm {SO}}(2)$ acts
on the set of side directions $(u_1, \dots, u_n)$ diagonally.

The polygon space $M_\ell$ emerges in topological robotics as the
configuration space of the planar linkage, a simple mechanism
consisting of $n$ bars of length $l_1, \dots, l_n$ connected by
revolving joints forming a closed planar polygonal chain. The
positions of two adjacent vertices are fixed but the other
vertices are free to move and the angles between the bars are
allowed to change.
\begin{center}
\resizebox{7cm}{4cm}
{\includegraphics[97,375][497,613]{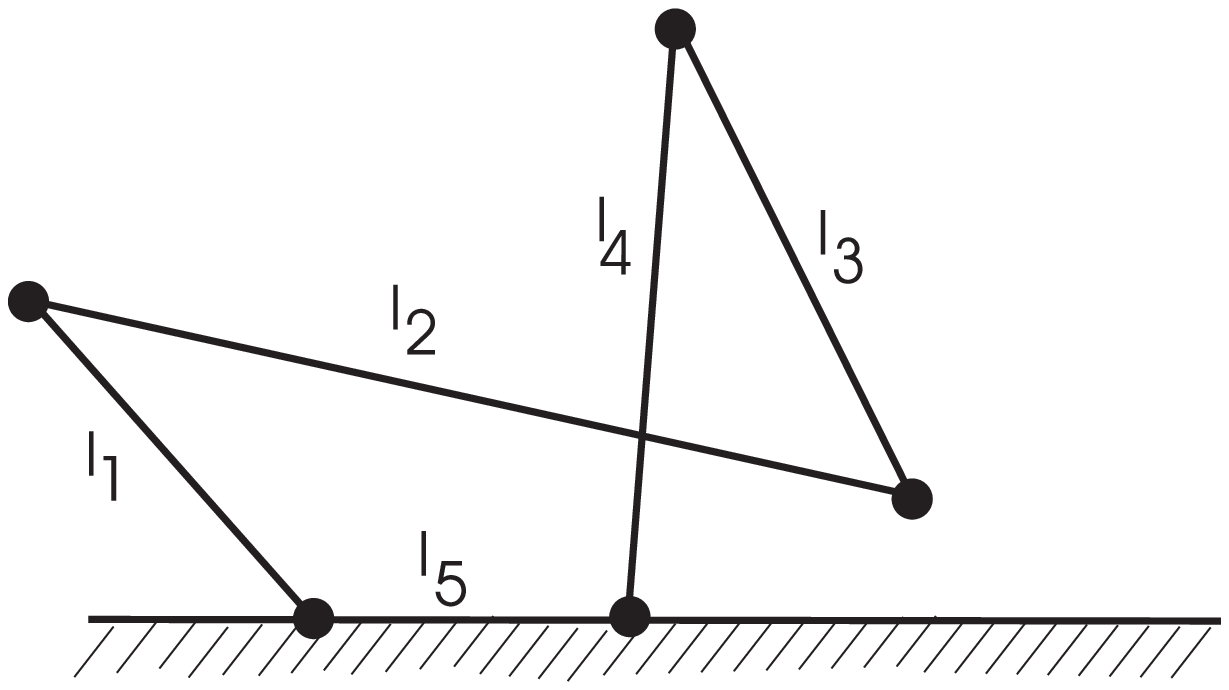}}
\end{center}
The spaces $M_\ell$ also appear in molecular biology where they
describe the space of shapes of closed molecular chains.

Statistical shape theory is another subject where the spaces
$M_\ell$ play a role: they describe the space of shapes having
certain geometric properties with respect to the central point,
see \cite{Kendall}.

The configuration space $M_\ell$ depends on the length vector
\begin{eqnarray}
\ell=(l_1, \dots,l_n)\in \R^n_+ \end{eqnarray} in an essential
way. Here $\R^n_+$ denotes the set of vectors in $\R^n$ having
nonnegative coordinates. Clearly, $M_\ell = M_{t\ell}$ for any
$t>0.$

The length vector $\ell$ is called {\it generic} if
$\sum\limits_{i=1}^n l_i\epsilon_i\not=0$ for any choice
$\epsilon_i=\pm 1$. It is known that for a generic length vector
$\ell$ the space $M_\ell$ is a closed smooth manifold of dimension
$n-3$. If the length vector $\ell$ is not generic then $M_\ell$ is
a compact $(n-3)$-dimensional manifold having finitely many
singular points.

The moduli spaces $M_\ell$ of planar polygonal linkages were
studied extensively by W. Thurston and J. Weeks \cite{TW}, K.
Walker \cite{Wa},
 A. A. Klyachko \cite{Kl},
M. Kapovich and J. Millson \cite{KM1}, J.-Cl. Hausmann and A.
Knutson \cite{HK} and others. The Betti numbers of $M_\ell$ as
functions of the length vector $\ell$ are described in \cite{FS};
we recall the result of \cite{FS} later in \S \ref{sec1}. A. A.
Klyachko \cite{Kl} found the Betti numbers of spacial polygon
spaces.

In this paper we view the length vector $\ell\in \R^n_+$ as a
random variable whose statistical behavior is described by
 a probability
measure $\mu$ on $\R^n_+$. The $p$-dimensional Betti number
$b_p(M_\ell)$ is then a function
\begin{eqnarray}
b_p: \R^n_+\to \Z, \quad  \ell \, \mapsto \, b_p(M_\ell)\in
\Z.\end{eqnarray}  Our goal is to find the asymptotics for $n$
large of the {\it average Betti numbers}, defined by
\begin{eqnarray}\label{barbi} b_p(n,\mu) =
\int_{\R^n_+} b_p(M_\ell)d\mu.\end{eqnarray} The choice of measure
$\mu$ reflects apriori information about the problem. We study in
detail two special choices for $\mu$: {\it
\begin{enumerate}
\item [(a)]  $\mu=\mu_a$ is the probability measure on $\R^n_+$
supported on the unit simplex $\Delta^{n-1}$ such that
$\mu|_{\Delta^{n-1}}$ coincides with the Lebesgue measure on
$\Delta^{n-1}\subset \R^n_+$ normalized so that
$\mu(\Delta^{n-1})=1$. Recall that $\Delta^{n-1}\subset \R^n_+$ is
described by the inequalities $l_i\geq 0$, $\sum l_i=1$.

\item[(b)] $\mu=\mu_b$ is the probability measure on $\R^n_+$
supported on the unit cube $\square^n\subset \R^n_+$ such that
$\mu|_{\square^n}$ is the Lebesgue measure. Here $\square^n$ is
given by the inequalities $0\leq l_i\leq 1$, $i=1, \dots, n$ and
$\mu(\square^n)=1$.
\end{enumerate}}

The main result of this paper states that in cases (a) and (b) for
any fixed $p \ge 0$ and large $n$ the following asymptotic formula
holds
\begin{eqnarray}\label{asym}
b_p(n, \mu) \sim \left(\begin{array}{c} n-1\\ p\end{array}
\right).
\end{eqnarray}
This can be expressed by saying that for large $n$ the random
manifold $M_\ell$ is connected, its first Betti number is $n-1$,
the second is $(n-1)(n-2)/2$ and so on. More precisely, our main
result is the following:

\begin{Theorem}\label{thm1} Let $\mu$ be either of the measures $\mu_a$ or $\mu_b$
described above. Then for any $p \ge 0$ there exist constants $C
>0$ and $0<a<1$ such that for any $n$ the average Betti number
$b_p(n,\mu)$, given by (\ref{barbi}), satisfies
\begin{eqnarray}\label{ineqthm1}
\left| b_p(n,\mu) - \left(\begin{array}{c} n-1\\ p\end{array}
\right)\right| \, < Ca^n.
\end{eqnarray}
\end{Theorem}

This result appears quite surprising for two reasons. Firstly, it
states that the asymptotic values of the average Betti number
$b_p(n, \mu)$ as $n\to \infty$ are equal for the measures $\mu_a$
and $\mu_b$ and raises the intriguing question about the
universality of the obtained asymptotic values.

Secondly, the binomial coefficient which appears in estimate
(\ref{ineqthm1}) equals the $p$-dimensional Betti number
$b_p(M_{\ell^{*}})$ of the configuration space of the equilateral
linkage $\ell^{*} = (1, \dots , 1)$ for any $p$ satisfying
$2p<n-3$, see \cite{FS}, Examples 3 \& 4. Hence the Betti numbers
of the moduli space $M_{\ell^{*}}$ approximate the average Betti
numbers $b_p(n, \mu)$ for $n$ large. On the other hand, by Theorem
2 of \cite{FS} the total Betti number, $\sum_{p=0}^{n}
b_p(M_{\ell}),$ viewed as a function of the length vector $\ell\in
\R^n_+$, is maximal for the equilateral linkage $\ell=\ell^{*}$.
This comparison shows that at least for some values of $p$ the
maximum of the individual Betti numbers $b_p(M_\ell)$ viewed as a
function of $\ell$ must be higher than $b_p(M_{\ell^\ast})$.

Indeed, this is the case since for the length vector
$\ell_\epsilon=(1,1,1, \epsilon, \dots, \epsilon)$, where
$\epsilon$ is a small positive number which appears $n-3$ times,
the Betti number $b_p(M_{\ell_\epsilon})$ equals
$$2\cdot \left(\begin{array}{c}n-3\\ p\end{array}
\right) = 2\cdot \frac{(n-2-p)(n-1-p)}{(n-1)(n-2)}\cdot
\left(\begin{array}{c}
n-1\\ p\end{array}\right)
$$
(see \cite{FS}, Example 2) which for large $n$ is nearly twice the
average Betti number $b_p(n,\mu)$. However if $p\sim n/2$
then $$b_p(M_{\ell_\epsilon})\sim \frac{1}{2}\left(\begin{array}{c} n-1\\
p\end{array}\right).$$

In a subsequent work we shall describe a generalization of Theorem
\ref{thm1} which allows the dimension $p$ to grow with $n$.

\section{Reduction of the problem}\label{sec1}

In this section we express the average Betti numbers in terms of
certain volumes.

To state this result we need to introduce some more notations. For
a subset $J\subset \{1, \dots, n\}$ we denote by $\phi_J: \R^n \to
\R$ the linear functional given by
\begin{eqnarray}\label{functional}
\phi_J(l_1, \dots, l_n) = \sum_{i\in J} l_i - \sum_{i\notin J} l_i
\end{eqnarray} and by $H_J$ the half-space
\begin{eqnarray}\label{hj}
H_J = \{\ell\in \R^n; \phi_J(\ell)<0\}.\end{eqnarray} Further,
$C^n$ denotes the cone
\begin{eqnarray}
C^n \, =\, \{(l_1, \dots, l_n); \, l_1\geq l_2\geq \dots l_n\geq
0\}\subset \R^n_+.\end{eqnarray}

The main result of this section provides a reduction of the
problem of computing the average Betti numbers.

\begin{Proposition}\label{prop2} Let $\mu$ be a probability measure on $\R^n_+$
having the following two properties: \begin{itemize}\item[(I)]
$\mu$ is invariant with respect to the action of the symmetric
group $\Sigma_n$ on $\R^n_+$ permuting coordinates; \item[(II)]
$\mu(L\cap \R^n_+)=0$ for any proper linear subspace $L\subset
\R^n$. \end{itemize} Then the average Betti number $b_p(n,\mu)$
equals
\begin{eqnarray}\label{propstatement}
b_p(n,\mu) \, =\, n!\cdot \sum_J \mu(H_J\cap C^n),
\end{eqnarray}
where $J$ runs over all subsets $J\subset \{1, 2, \dots, n\}$
containing $1$ and having cardinality either $|J|=p+1$ or
$|J|=n-2-p$.
\end{Proposition}

\begin{proof}
First, recall the result of \cite{FS} stating that the Betti
numbers $b_p(M_\ell)$, as functions of the length vector
$\ell=(l_1, \dots, l_n)$, can be computed by counting certain
subsets of the index set $\{1, \dots, n\}$. A subset $J\subset
\{1, \dots, n\}$ is called {\it short} if
\begin{eqnarray*}
\sum_{i\in J}l_i \, < \, \sum_{i\notin J}l_i.\end{eqnarray*} A
subset is called {\it long} if its complement is short. A subset
$J\subset \{1, \dots, n\}$ is called {\it median} if
\begin{eqnarray*}
\sum_{i\in J} l_i=\sum_{i\notin J} l_i.\end{eqnarray*}

Fix an index $1\leq i\leq n$ such that $l_i$ is maximal among
$l_1, \dots, l_n$. Denote by $a_p(\ell)$ the number of short
subsets $J\subset \{1, \dots, n\}$ of cardinality $|J|=1+p$
containing $i$. Denote by $\tilde a_p(\ell)$ the number of median
subsets $J\subset \{1, \dots, n\}$ containing $i$ and such that
$|J|=1+p$. It was proven in \cite{FS} that for $p=0, 1, \dots,
n-3$ one has
\begin{eqnarray}\label{fs}
b_p(M_\ell) = a_p(\ell) + \tilde a_p(\ell) + a_{n-3-p}(\ell).
\end{eqnarray}

It is easy to see that the manifolds $M_{\ell_1}$ and $M_{\ell_2}$
are diffeomorphic if the length vector $\ell_1\in \R^n_+$ is
obtained from $\ell_2\in \R^n_+$ by permuting the components. In
other words, the order of the coordinates $l_1, \dots, l_n$ in the
length vector $\ell=(l_1, \dots, l_n)$ is irrelevant.

Let $\mu$ be a probability measure on $\R^n_+$ having properties
(I) and (II). Property (II) implies that
\begin{eqnarray}
\int_{\R^n_+} \tilde a_p(\ell)d\mu \, =\, 0
\end{eqnarray}
since the function $\ell\mapsto\tilde a_p(\ell)$ is zero on the
complement of a union of finitely many linear hyperplanes.
Integrating (\ref{fs}) we have
\begin{eqnarray}\label{bpmu}
b_p(n, \mu)= a_p(n, \mu) + a_{n-3-p}(n,\mu).
\end{eqnarray}
where
\begin{eqnarray}\label{apmu}
 a_p(n, \mu) \, =\,  \int_{\R^n_+} a_p(\ell)d\mu \, =\,
 n!\cdot \int_{C^n} a_p(\ell) d\mu.
 \end{eqnarray}
The second equality follows from property (I) of $\mu$. The
function $a_p|_{C^n}$ is quite simple. Denote by $\sigma$ the step
function
\begin{eqnarray*}
\sigma(x) = \left\{
\begin{array}{lll}
1, &\mbox{if} & x<0,\\
0, &\mbox{if} & x\geq 0,
\end{array}\right.
\end{eqnarray*}
where $x\in \R$. Then we may write
\begin{eqnarray}
a_p(\ell) = \sum_J \sigma((\phi_J(\ell)))
\end{eqnarray}
where $J$ runs over all subsets of $\{1, \dots, n\}$ of
cardinality $p+1$ containing $1$. Integrating we find
\begin{eqnarray}\label{apell}
a_p(n,\mu) = n! \cdot \int_{C^n} a_p(\ell)d\mu  \nonumber\\
 = n! \cdot \sum_J
\int_{C^n} \sigma(\phi_J(\ell))d\mu.
\end{eqnarray}
In formula (\ref{apell}) $J$ runs over all subsets $J\subset \{1,
\dots, n\}$ satisfying $1\in J$, $|J|=p+1$. Obviously
\begin{eqnarray}\label{vj}
 \int_{C^n}\sigma(\phi_J(\ell))d\mu \, =\, \mu(H_J\cap
C^n)
\end{eqnarray}
where $H_J$ is given by (\ref{hj}). Our statement
(\ref{propstatement}) now follows by combining (\ref{bpmu}),
(\ref{apmu}), (\ref{apell}) and (\ref{vj}).
\end{proof}

\section{Simplices and volumes}

In this section we give a geometric interpretation of the
quantities appearing on the RHS of formula (\ref{propstatement}).

Recall that we denote by $\mu_a$ the measure on $\R^n_+$ with
support on the unit simplex $\Delta^{n-1}\subset \R^n_+$ such that
$\mu_a|_{\Delta^{n-1}}$ coincides with the normalized Lebesgue
measure, $\mu_a(\Delta^{n-1}) =1$. Similarly, $\mu_b$ denotes the
probability measure on $\R^n_+$ supported on the unit cube
$\square^n\subset \R^n_+$ such that $\mu_b|_{\square^n}$ coincides
with the Lebesgue measure.

Denote by $A\subset \R^n_+$ the simplex of dimension $n$ having
the vertices
\begin{eqnarray}
c_0&=&(0, \dots, 0),\nonumber\\
c_1&=&(1, 0, \dots, 0),\nonumber\\
c_2&=& \frac{1}{2}(1,1,0, \dots, 0),\label{verticesA}\\
\dots\nonumber\\
c_n &=& \frac{1}{n}(1,1,\dots, 1).\nonumber
\end{eqnarray}
Similarly, denote by $B\subset \R^n_+$ the simplex of dimension
$n$ having the vertices
\begin{eqnarray}\label{verticesB}
c'_i =(\underbrace{1, \dots, 1}_{i \, \, \mbox{times}}, 0\dots,0)=
i\cdot c_i,
\end{eqnarray}
where $i=0, 1, \dots, n$.
\begin{Proposition}\label{prop3} For any subset $J\subset \{1, \dots, n\}$ one has
\begin{eqnarray}
n!\cdot \mu_a(H_J\cap C^n) = \frac{\vol(H_J\cap A)}{\vol(A)},
\quad n!\cdot \mu_b(H_J\cap C^n) = \frac{\vol(H_J\cap
B)}{\vol(B)}.
\end{eqnarray}
In other words, the quantities appearing in Proposition
\ref{prop2} can be interpreted as ratios of volumes of certain
simplices and their parts cut off by a half-space.
\end{Proposition}
\begin{proof} One checks that the intersection $\Delta^{n-1}\cap C^n$ coincides with
the $(n-1)$-dimensional simplex $A'\subset A$ with vertices $c_1,
\dots, c_n$ and the intersection $\square^n\cap C^n$ equals $B$.
Clearly
$$\mu_a(A') = \frac{\vol_{n-1}(A')}{\vol_{n-1}(A)}=(n!)^{-1}, \quad \mu_b(B)= \vol(B) =
(n!)^{-1}.$$ Here $\vol_{n-1}$ denotes the $(n-1)$-dimensional
Euclidean area. Hence we find that
\begin{eqnarray}
n!\cdot \mu_a(H_J\cap C^n) = \frac{\vol_{n-1}(H_J\cap
A')}{\vol_{n-1}(A')} = \frac{\vol(H_J\cap A)}{\vol(A)}
\end{eqnarray}
and
\begin{eqnarray}
n!\cdot \mu_b(H_J\cap C^n) = \frac{\vol(H_J\cap B)}{\vol(B)}.
\end{eqnarray}
This completes the proof.
\end{proof}

We are led to consider the following simple geometric problem.

\begin{Proposition}\label{cutvol}
Given $n+1$ points $v_0, \dots, v_n\in \R^n$ in general position
and a linear functional $\phi: \R^n\to \R$, consider the simplex
$\Sigma$ spanned by $v_0, \dots, v_n$ and its intersection with
the half space $H=\{v\in \R^n; \phi(v)<0\}$. Denote by
\begin{eqnarray}
q_i=\phi(v_i), \quad i=0, \dots, n
\end{eqnarray}
the values of the functional $\phi$ on the vertices of $\Sigma$.
Assume that the numbers $q_0, q_1, \dots, q_n$ are all distinct
and the vertices are labelled so that $q_i<0$ for all $i=0,\dots,
m$ and $q_i\geq 0$ for $i=m+1, \dots, n$. Then the ratio of the
volumes
\begin{eqnarray}\label{ratio2} r= \frac{{\vol}(\Sigma\cap
H)}{{\vol}(\Sigma)}\end{eqnarray} equals
\begin{eqnarray}\label{slice}
r = \sum_{i=0}^m \, \prod\limits_{\begin{array}{c}0\leq j\leq n\\
j\not= i\end{array}}\, \frac{q_i}{q_i- q_j}.
\end{eqnarray}
\end{Proposition}

Proposition \ref{cutvol} is well known. Different expressions for
the volume cut off a simplex by a half space were obtained in
\cite{CS} (see Theorem 2 of \cite{CS}) and more recently in
\cite{A}, \cite{Ge}, \cite{Va}. Proposition \ref{cutvol} can be
easily obtained from the results mentioned above. For convenience
of the reader we briefly sketch the proof.

\begin{proof}[Proof of Proposition \ref{cutvol}]  Consider for
$x\in \R$ the half space $H_x=\{v\in \R^n; \phi(v)<x\}$ and the
real valued function
\begin{eqnarray}
r(x) = \frac{{\vol}(H_x\cap B )}{{\vol}(B)}.
\end{eqnarray}
Without loss of generality we will assume that $q_0<q_1<
\dots<q_n$; for convenience we set $q_{-1}=-\infty$ and
$q_{n+1}=+\infty$.

The function $r(x)$ has the following properties:
\begin{enumerate}
\item[(i)] when restricted on each subinterval $[q_{i-1}, q_i]$,
the function $r(x)$ is a polynomial $p_i(x)$ of degree $n$ with
real coefficients where $i=0, \dots, n+1$.

\item[(ii)] $p_0(x)\equiv 0$ and $p_{n+1}(x)\equiv 1$.

\item[(iii)] $r(x)$ has continuous derivatives up to order $ n-1$.
\end{enumerate}

Properties (i) and (iii) (proven in \cite{CS}) imply that for any
$0 \le i \le n$ there exists $\beta_{i} \in \R$ such that
\begin{eqnarray}
p_{i+1}(x) - p_i(x) = \beta_i(x-q_i)^n.
\end{eqnarray}
Hence
\begin{eqnarray}\label{pii}
p_{i+1}(x) = \sum_{k=0}^i \beta_k(x-q_k)^n
\end{eqnarray}
where the coefficients $\beta_i$ can be found from the polynomial
identity
\begin{eqnarray}\label{system}
\sum_{k=0}^n \beta_k(x-q_k)^n\equiv 1.
\end{eqnarray}
Equation (\ref{system}) has a unique solution given by
\begin{eqnarray}\label{alphai}
\beta_k = (-1)^n \cdot\prod_{\begin{array}{c} j=0\\
j\not=k\end{array}}^n (q_k-q_j)^{-1}, \quad k=0, 1, \dots, n.
\end{eqnarray}
Indeed, by comparison of coefficients in (\ref{system}) we obtain
the linear system
\begin{eqnarray}
\sum_{k=0}^n \beta_k q_k^{n-i} = \left\{
\begin{array}{lll}
(-1)^n & \mbox{if} & i=0,\\  0 & \mbox{if} & i=1, \dots, n.
\end{array} \right.
\end{eqnarray}
which can be written in the matrix form
\begin{eqnarray}
\left[
\begin{array}{cccc}
1 & 1& \dots & 1\\
q_0 & q_1& \dots & q_n\\
\dots& \dots&\dots &\dots\\
q_0^{n-1}& q_1^{n-1} & \dots & q_n^{n-1}\\
q_0^n& q_1^n & \dots & q_n^n
\end{array}
\right] \left[
\begin{array}{c}
\beta_0\\
\beta_1\\
\dots\\
\beta_{n-1}\\
\beta_n\end{array} \right] = \left[\begin{array}{c} 0\\ 0\\
\dots\\0\\ (-1)^n\end{array} \right]
\end{eqnarray}
To solve it one applies Cramer's rule. The determinant of the
system is the Vandermonde determinant $\prod_{j>l} (q_j-q_l)$ and
the numerator of the fraction expressing $\beta_k$ can be computed
to be the Vandermonde determinant $$(-1)^k\, \cdot
\prod_{\begin{array}{c} j>l\\ j, l\not= k\end{array}} (q_j-q_l)$$
so that we obtain
$$\beta_k= (-1)^k \prod_{k>l}(q_k - q_l)^{-1}\cdot
\prod_{j>k}(q_j-q_k)^{-1} = (-1)^{n}\prod_{j \not=
k}(q_k-q_j)^{-1}$$ which coincides with (\ref{alphai}).
Substituting the obtained value into (\ref{pii}) we find that for
$x\in [q_{i-1},q_i]$ one has
\begin{eqnarray}\label{slice1}
r(x) \, = \, p_i(x) = \sum_{k=0}^{i-1}
\prod\limits_{\begin{array}{c} 0\leq j\leq n\\
j\not=k\end{array}} \, \frac{q_k-x}{q_k-q_j}.
\end{eqnarray}
Formula (\ref{slice}) is obtained from (\ref{slice1}) by setting
$x=0$.
\end{proof}

For our applications we need to have a more general formula for
the ratio (\ref{ratio2}) covering the case when some values
$q_i=\phi(c_i)$ coincide. Formula (\ref{slice}) is not well
defined in this case since some of the denominators might vanish.
Our setting is as follows. Let $v_0, \dots, v_n\in \R^n$ be points
in general position spanning a simplex $B$. Let $\phi: \R^n\to \R$
be a linear functional. We consider the ratio of the volumes
(\ref{slice}) where $H$ is the half-space $H=\{v\in \R^n; \phi(v)
< 0\}$. Denote $q_i=\phi(v_i)\in \R$, where $i=0, \dots, n$.
Suppose that there is a decomposition of the set of indices into
disjoint subsets
$$\{0, 1, \dots, n\} = \bigsqcup_{l=0}^s I_l$$
such that $q_i=Q_l$ for all $i\in I_l$ and the numbers $Q_0,
\dots, Q_s\in \R$ are pairwise distinct. We denote by $k_l$ the
number $|I_l|-1$. Thus, the multiplicity of $Q_l$ is $k_l+1$,
where $l=0,1, \dots, s$ and one has
\begin{eqnarray}\label{sum5}
k_0+k_1+\dots+k_s =n-s.
\end{eqnarray}
Given two nonnegative integers $s$ and $a$, we denote by $P(s,a)$
the set of all functions $$\delta: \{0,1, \dots, s\}\to \Z_{\ge
0}, \quad i\mapsto \delta_i$$ satisfying $\sum_{j=0}^s\delta_j
=a.$ The set $P(s,a)$ labels partitions of $a$ in $s+1$ summands.

\begin{Proposition}\label{ratio3}
Assume that $Q_i<0$ for all $i=0,\dots, m$ and $Q_i\geq 0$ for all
$i=m+1, \dots, s$. Then the ratio of volumes (\ref{ratio2}) equals
\begin{eqnarray}\label{ratiogeneral}
r= \sum_{i=0}^m \,\left[ F_i\cdot \prod\limits_{\begin{array}{c} 0\leq j\leq s\\
j\not= i\end{array}}\left(
\frac{Q_i}{Q_i-Q_j}\right)^{k_j+1}\right]
\end{eqnarray}
where $F_i$ equals
\begin{eqnarray}\label{fi}
\sum\limits_{\delta\in P(s, k_i)}
\left(
\begin{array}{c} n\\ \delta_i\end{array}
\right)\cdot (-Q_i)^{k_i-\delta_i}\cdot
\prod\limits_{\begin{array}{c} 0\leq j\leq s\\ j\not=i\end{array}}
\left(\begin{array}{c} k_j+\delta_j\\
\delta_j\end{array}\right)\cdot (Q_i-Q_j)^{-\delta_j}.
\end{eqnarray}
\end{Proposition}

Note that  if $k_i=0$, then $P(s,k_{i})= \{ \delta \equiv 0\}$ and
hence $F_i=1$. In particular, if $s = n$ and therefore $k_{i} = 0$
for any $0 \le i \le n$ formula (\ref{fi}) coincides with formula
(\ref{slice}). Moreover, we point out that the coefficients
$F_{i}$ are homogenous of degree $0$ in the variables $Q_{0},
\dots , Q_{s}.$

\begin{proof}[Proof of Proposition \ref{ratio3}]
Let $\psi: \R^n \to \R$ be a  linear functional, a small
perturbation of $\phi: \R^n\to \R$, such that the values
$p_i=\psi(v_i)$ are all distinct. Without loss of generality we
may assume that $p_i<0$ for $0\leq i \leq m'$ and $p_i>0$ for $m'
<i \leq n$ where $m'=m +1+\sum_{i=0}^m k_i$. Denote $H'=\{v\in
\R^n; \psi(v)<0\}$. The perturbed ratio
\begin{eqnarray}\label{ratioperturb}
r'= \frac{\vol(H'\cap \Sigma)}{\vol(\Sigma)}
\end{eqnarray}
tends to $r$ when $\psi$ tends to $\phi$. By Proposition
\ref{cutvol} it can be written in the form
\begin{eqnarray}\label{ratioper}
r' = \sum_{i=0}^{m'} \left[ p_i^n \cdot \prod\limits_{\begin{array}{c} j\not= i\\
0\leq j \leq n\end{array}}\left(p_i - p_j\right)^{-1}\right].
\end{eqnarray}
To find its limit we will use the following general formula, see
\cite{MT},
\begin{eqnarray}\label{formula}
\sum_{i=0}^k f(x_i)\cdot \prod_{\begin{array}{c} 0\leq j\leq k\\
j\not= i\end{array}} (x_i-x_j)^{-1} = \frac{1}{k!} f^{(k)}(\xi)
\end{eqnarray}
where $f(x)$ is a real valued smooth function, $x_0, \dots, x_k$
are distinct real numbers and $\xi$ is a number lying in the
smallest interval containing these points. Applying
(\ref{formula}) $(m+1)$ times to (\ref{ratioper}) and passing to
the limit as $\{p_0, \dots, p_n\} \to \{Q_0, \dots, Q_s\}$ we
obtain
\begin{eqnarray}\label{multiple}
r = \sum_{i=0}^m \frac{1}{k_i!}\left[x^n\cdot
\prod_{\begin{array}{c} 0\leq l\leq s\\ l\not=i\end{array}}
(x-Q_l)^{-k_l-1}\right]^{(k_i)}_{x=Q_i}.
\end{eqnarray}
Recall the following formula for higher derivatives of products
$$(f_0f_1\cdots f_s)^{(k)} = \sum_{\delta\in P(s,k)}
\frac{k!}{\delta_0!\cdots \delta_s!}\, \,  f_0^{(\delta_0)}\cdots
f_s^{(\delta_s)}.$$ Applying it to (\ref{multiple}) and using that
$n=k_i+\sum_{j\not= i} (k_j+1)$ one obtains (\ref{ratiogeneral})
after certain elementary transformations.
\end{proof}

\section{Sequences of densities}

By Proposition \ref{prop2} combined with Proposition \ref{prop3},
to compute the average Betti numbers $b_p(n,\mu)$ one has to know
the volumes cut off a simplex by certain half-spaces. The result
of Propositions \ref{cutvol} and \ref{ratio3} show that to find
these volumes it is enough to know the values of the functionals
determining the half-spaces on the vertices of the simplices $A$
and $B$. We investigate these values in this section.

Let $J\subset \{1, \dots, n\}$ be an arbitrary subset and $\phi_J:
\R^n\to \R$ be the linear functional (\ref{functional}). The
values of the functional $\phi_J$ on the vertices $c_i$ of the
simplex $A$ (see (\ref{verticesA})) equal $\phi_J(c_0) =0$ and for
$i\geq 1$
\begin{eqnarray}\label{phi}
\phi_J(c_i) = \frac{1}{i} \left|J\cap \{1, \dots, i\}\right|\,-\,
\frac{1}{i}\left|\bar J\cap \{1, \dots, i\}\right|
\end{eqnarray}
where $\bar J$ denotes the complement of $J$ in $\{1, \dots ,
n\}$. Let
\begin{eqnarray}\label{dens}
\alpha_i(J) = \frac{1}{i}\left|J\cap \{1, \dots, i\}\right|, \quad
i=1, \dots, n
\end{eqnarray}
denote  the {\it density} of the set $J$ in the interval $\{1,
\dots, i\}$. Clearly $0 \, \leq \, \alpha_i(J)\, \leq\,  1$.

\begin{lemma}\label{lm1} Let $J\subset \{1, \dots, n\}$ be a subset of
cardinality $|J|=p\geq 1$. Then the following estimates hold:
\begin{enumerate}
\item[(a)]  For $2p \le i \le n$ one has
\begin{eqnarray}\label{half}
\alpha_i(J) \leq \frac{1}{2} \end{eqnarray} and equality in
(\ref{half}) may only hold for $i=2p$. \item[(b)] If $1\leq i,
j\leq 2p$, then either $\alpha_i(J)=\alpha_j(J)$ or
$$\left|\alpha_i(J) - \alpha_j(J)\right| \geq \frac{1}{(2p)^2}.$$
\item[(c)] If  $\alpha_i(J) <\frac{1}{2}$ for some $1 \le i \le
n$, then
$$\alpha_i(J) \leq \frac{1}{2} - \frac{1}{2(2p+1)}.$$
\item[(d)] For any $8p^3 \le i\leq n$ one has
\begin{eqnarray}
0\leq \alpha_i(J) \leq \frac{1}{8p^2}.
\end{eqnarray}
\end{enumerate}
\end{lemma}
\begin{proof}
(a) Let $k_i$ denote $|J\cap \{1, 2, \dots, i\}|$. Then for $i\geq
2p$ one has
$$\alpha_i(J) = \frac{k_i}{i} \leq \frac{p}{i} \leq \frac{p}{2p} =1/2.$$
If $\alpha_i(J) = 1/2,$ then the above inequalities imply that $i = 2p.$

(b) Suppose now that  $1\leq i, j\leq 2p$. Then
\begin{eqnarray}\label{dif}
\left|\alpha_i(J)-\alpha_j(J)\right| = \left|\frac{k_i}{i} -
\frac{k_j}{j}\right|=
\frac{\left|jk_i-ik_j\right|}{ij}.\end{eqnarray} We obtain that
either (\ref{dif}) vanishes or $\left|jk_i - ik_j\right|\geq 1$
and hence (\ref{dif}) is greater or equal to
$$\frac{1}{ij}\geq \frac{1}{(2p)^2}$$
proving (b).

(c) Note that for $i\geq 2p+2$ one has
$$\alpha_i(J) =\frac{k_i}{i} \leq \frac{p}{2p+2} = \frac{1}{2} -
\frac{1}{2(p+1)}.$$
Now consider the case $i\leq 2p+1.$ By assumption $\alpha_i(J) <1/2$, i.e.
$2k_i<i$ or $2k_i\leq i-1.$ It implies that
$$\alpha_i(J) =\frac{k_i}{i} \leq \frac{i-1}{2i} = \frac{1}{2} -
\frac{1}{2i} \leq \frac{1}{2} - \frac{1}{2(2p+1)}.$$

(d) Note that if $i \ge 8p^3$ and $p \ge 1,$ then
$$\alpha_i(J) =\frac{k_i}{i}\leq
\frac{p}{i}\leq \frac{1}{8p^2}.$$
 \end{proof}

\begin{lemma}\label{lm3}
Let $J \subset \{1, \dots , n\}$ with $|J|=p\geq 1$. A nonzero
number may appear in the sequence of densities $\alpha_1(J),
\dots, \alpha_n(J)$
 at most $p$ times.
\end{lemma}

\begin{proof}
 The densities
$\alpha_i(J)$ satisfy the following recurrent relation:
\begin{eqnarray}
\alpha_{i+1}(J) = \left\{
\begin{array}{lll}
\frac{i}{i+1} \alpha_i(J),& \mbox{if} & i+1\notin
J,\\ \\
\frac{i}{i+1} \alpha_i(J) +\frac{1}{i+1},& \mbox{if} & i+1\in J.\\
\end{array}
\right.
\end{eqnarray}
It follows that $\alpha_{i+1}(J)<\alpha_i(J)$ if $\alpha_i(J)>0$
and $i+1\notin J$. On the other hand, $\alpha_{i+1}(J)\geq
\alpha_i(J)$ if $i+1\in J$. Hence, for $i<j$, the equality
$\alpha_i(J)=\alpha_j(J)>0$ implies that at least one of the
intermediate indices $i+1, \dots, j$ belongs to $J$.

Assume now that for $i_1<i_2<\dots < i_k$ one has $\alpha_{i_1}(J)
= \dots = \alpha_{i_{k}}(J)> 0$. The set $J$ divides $\{1, \dots,
n\}$ into $p+1$ subintervals, i.e. subsets of consecutive integers
in $\{1, \dots, n\}$ all of which are not in $J$ except the first
one. The leftmost interval contains no elements of $J$ and might
be empty. As explained above, each of the subintervals may contain
at most one of the integers $i_1, \dots, i_k$. The density
$\alpha_i(J)$ vanishes iff $i$ lies in the leftmost subinterval.
This shows that $k\leq p$ as claimed.
\end{proof}

The following examples show that Lemma \ref{lm3}
cannot be improved.

\begin{Example}{\rm
(i) Consider $J = \{ 1,2,\dots, p\}$. Then $\alpha_i(J)=1$ if and
only if $1 \le i \le p$, i.e. the multiplicity of the value $1$ is
$p$.

(ii) Let $J=\{1, 3, 5, \dots, 2p-1\}$. Then $\alpha_i(J)$ equals
$1/2$ exactly $p$ times.

(iii) Suppose that $J=\{n-p+1, n-p+2, \dots, n\}$. Then
$\alpha_i(J)=0$ for $i=1, 2, \dots, n-p$, i.e. the multiplicity of
the value zero is $n-p$. This shows that the bound of the
multiplicity of Lemma \ref{lm3} does not hold for the value zero.}
\end{Example}

The values \begin{eqnarray}\label{values}
q_{i}=\phi_J(c_i)\end{eqnarray} of the functional $\phi_J$ on the
vertices of the simplex $A$ (see (\ref{verticesA})) will play an
important role in the sequel. Clearly $q_0=0$ and
\begin{eqnarray} q_i = 2\alpha_i(J) -1\, \, \in [-1, 1], \quad i=1, \dots,
n,
\end{eqnarray}
see (\ref{phi}). Let us restate Lemmas \ref{lm1} and \ref{lm3} in
terms of the $q_i$'s:

\begin{lemma}\label{lm2} Let $J\subset \{1, \dots, n\}$ be a subset of
cardinality $|J|=p\geq 1$. Then:
\begin{enumerate}
\item[(a)]    For $2p \le i \le n$ one has $q_i\leq 0$ with
equality possible only for $i=2p$. \item[(b)] If $1\leq i, j\leq
2p$ then either $q_i=q_j$ or
$$\left|q_i - q_j\right| \geq \frac{1}{2p^2}.$$
\item[(c)] If $q_i < 0 $ for some $1 \le i \le n$  then
$$q_i \leq - \frac{1}{2p+1}.$$
\item[(d)] For $8p^3 \leq i\leq n$ one has
\begin{eqnarray}
-1\leq q_i \leq -1+ \frac{1}{4p^2}.
\end{eqnarray}
\item[(e)] A number distinct from $-1$ may appear in the sequence
$q_0, \dots, q_n$ at most $p$ times.
\end{enumerate}
\end{lemma}

Next we consider the values $q'_i=\phi_J(c'_i)$ of the functional
$\phi_J$ on the vertices of simplex $B$ (see (\ref{verticesB})).

\begin{lemma}\label{lmqprime}
Let $J\subset \{1, \dots, n\}$ be a subset of cardinality
$|J|=p\geq 1$. Then:
\begin{enumerate} \item[(a)] The numbers $q'_i$ satisfy $-i\le
q'_i\le i$; \item[(b)]    For $2p \le i \le n$ one has $q'_i\leq
0$ with equality possible only for $i=2p$. \item[(c)] If
$q'_i\not=q'_j$, then $\left|q'_i - q'_j\right| \geq 1.$
\item[(d)] If $q'_i < 0 $ for some $1 \le i \le n$,  then $q'_i
\leq - 1.$ \item[(e)] A number may appear in the sequence $q'_0,
\dots, q'_n$ at most $p+1$ times.
\end{enumerate}\end{lemma}
\begin{proof}
Since $q'_i=iq_i$, statements (a) - (d) follow from Lemma
\ref{lm2} and from the observation that $q'_i$ is an integer. To
prove (e) we note that the numbers $q'_i$ satisfy the following
recurrent relation:
$$q'_{i+1} = \left\{
\begin{array}{ll}
q'_i -1 & \mbox{if}\, \, i+1\notin J,\\ \\
q'_i +1 & \mbox{if}\, \,
i+1\in J.
\end{array}
\right.$$ Hence the sequence $q'_0, q'_1, \dots, q'_n$ has exactly
$p$ jumps up and decays between the jumps. This proves (e).
\end{proof}

\begin{Example}
{\rm Consider the sequence $J=\{1, 3, \dots, 2p-1\}$. Then (i)
$q'_i=0$ for $i=0, 2, 4, \dots, 2p$, (ii) $q'_i=1$ for $i=1, 3,
\dots 2p-1$ and (iii) $q'_i<0$ for $i>2p$. We see that in this
case zero appears in the sequence $q'_0, \dots, q'_n$ exactly
$p+1$ times.}
\end{Example}

\section{Proof of Theorem \ref{thm1} for $\mu=\mu_a$}

In this section we prove Theorem \ref{thm1} for the measure
$\mu=\mu_a$ described before the statement of Theorem \ref{thm1}.

Fix a subset $J\subset\{1, 2, \dots, n\}$ of cardinality $p\geq 1$
where we think of $p$ as being fixed and of $n$ as being large.
Let $A$ be the simplex with vertices (\ref{verticesA}). Our first
goal is to estimate the ratios of the form
\begin{eqnarray}
r_J= \frac{\vol(H_J\cap A)}{\vol (A)}
\end{eqnarray}
where $H_J\subset \R^n$ is the half-space $H_J=\{v\in \R^n;
\phi_J(v)<0\}$. The average Betti numbers are sums of ratios of
this kind, see Propositions \ref{prop2} and \ref{prop3}. Recall
that $\phi_J$ denotes the linear functional $\phi_J: \R^n\to \R$
given by (\ref{functional}).

The values of the functional $\phi_J$ on the vertices $c_i$ of $A$
are described in Lemma \ref{lm2}. In particular, by statement (a)
of Lemma \ref{lm2}, for large $n$ the majority of the values $q_i$
are negative, more precisely, at most $2p$ of them are positive.
Hence, one may expect that the volume of the section $H_J\cap A$
approximates $\vol(A)$ for large $n$. This is indeed the case.

To estimate the difference $1-r_J$ from above we will apply
Propositions \ref{cutvol} and \ref{ratio3}; our aim is to show
that it is exponentially small.

\begin{Proposition}\label{estim2} Given an integer $p\geq 1$,
there exist constants $C>0$ and
$0<a<1$ such that for any $n\geq p$ and any subset $J\subset \{1,
2, \dots, n\}$ of cardinality $p$ one has
\begin{eqnarray}\label{estim1}
1-r_J < C\cdot a^n.
\end{eqnarray}
\end{Proposition}

It will be apparent from the proof that one can take for $a$ an
arbitrary number satisfying $(2p+1)/(2p+2)<a<1$.

\begin{proof}[Proof of Proposition \ref{estim2}]
Consider the values $q_i=\phi_J(c_i)$, where $i=0,
1, \dots, n$. They may have multiplicities, i.e. the same value
may appear several times. We will denote by $Q_0> Q_1 >, \dots,
> Q_s\in [-1, 1]$ the different values of the sequence $q_0, \dots, q_n$.
Then there is a surjective mapping $\tau: \{0,1, \dots, n\}\to
\{0,1,\dots, s\}$ such that $q_i=Q_{\tau(i)}$. For $0\leq i\leq s$
we denote by $k_i+1$ the cardinality of the preimage
$\tau^{-1}(i)$. Clearly,
$$\sum_{i=0}^s k_i = n-s.$$
If $Q_i\leq 0$ for all $0\leq i\leq s$ then $H_J\cap A= A$, hence
$r_J=1$ and therefore (\ref{estim1}) trivially holds. Thus without
loss of generality we may assume that for some $0\leq m\leq s$ one
has $Q_i>0$ for $i=0, \dots, m$ and $Q_i\leq 0$ for $i= m+1,
\dots, s$. By statement (a) of Lemma \ref{lm2} we have
\begin{eqnarray}\label{mki}
\sum_{i=0}^m (k_i+1)\, \leq 2p.
\end{eqnarray}
We also have
\begin{eqnarray}\label{ineq3}
k_i\leq p, \quad i=0, 1, \dots, s-1,\end{eqnarray} see (e) of
Lemma \ref{lm2}. Inequality (\ref{ineq3}) also applies to the
multiplicity $k_s$ if $Q_s$ is distinct from $-1$.

Let $b=1+(2p)^{-2}.$ From statements (a), (b), and (d) of Lemma
\ref{lm2} we obtain that for any $0\leq i\leq m$ and $j\not=i$ one
has
\begin{eqnarray}\label{lowerbound}
|Q_i-Q_j|\, \geq\,  \left\{
\begin{array}{ll}
(2p^2+1)^{-1},& \mbox{for}\, \,0 \leq  j\leq 8p^3,\\ \\
(2p^2)^{-1} - (-1 +\frac{1}{4p^2})=b, & \mbox{for}\, \, 8p^3 <
j\leq n.
\end{array}
\right.
\end{eqnarray}
Let us explain the first line of this inequality. If $j$ satisfies
$0\leq j\leq 2p$ then the statement follows from (a) of Lemma
\ref{lm2}. If $j>2p$ then $Q_j<0$ and hence $Q_i-Q_j\geq
(2p+1)^{-1}\geq (2p^2+1)^{-1}$ by (c) of Lemma \ref{lm2}.

By applying Proposition \ref{ratio3} to $-\phi_{J}$ we have

\begin{eqnarray}\label{vol-formula}
1-r_J= \sum_{i=0}^m \, \, \left[F_i\cdot\prod\limits_{\begin{array}{c} 0\leq j\leq s\\
j\not= i\end{array}} \left(
\frac{Q_i}{Q_i-Q_j}\right)^{k_j+1}\right],
\end{eqnarray}
where $F_i$ is given by
\begin{eqnarray}\label{fi1}
\sum\limits_{\delta\in P(s, k_i)}\left(
\begin{array}{c} n\\ \delta_i\end{array}
\right)\cdot(-Q_i)^{k_i-\delta_i}\cdot
\prod\limits_{\begin{array}{c} 0\leq j\leq s\\ j\not=i\end{array}}
\left(\begin{array}{c} k_j+\delta_j\\
\delta_j\end{array}\right)\cdot (Q_i-Q_j)^{-\delta_j}.
\end{eqnarray}

We claim that for any $0 \le i \le m,$
\begin{eqnarray}\label{estim}
|F_i| \leq  (2p^2+1)^p n^{2p}.
\end{eqnarray}
Indeed, observe that by (\ref{lowerbound}),
$$\left| Q_i-Q_j \right| ^{-\delta_j} \le (2p^{2}+1)^{\delta_{j}}.
$$
To estimate the
binomial coefficient $\left(\begin{array}{c} k_j+\delta_j\\
\delta_j\end{array}\right)$ note that $\delta_j\le k_i$ and thus
$k_j+\delta_j\leq k_j+k_i\leq n$ and therefore
\begin{eqnarray}\label{binom1}
\left(\begin{array}{c} k_j+\delta_j\\
\delta_j\end{array}\right) \le n^{\delta_{j}}\quad \mbox{and}\quad
\left(\begin{array}{c} n\\
\delta_i\end{array}\right) \le n^{\delta_{i}} .
\end{eqnarray}
As $|Q_{i} | \le 1,$ each term in the sum (\ref{fi1}) can be
estimated as follows
\begin{eqnarray*}
\left(
\begin{array}{c} n\\ \delta_i\end{array}
\right)\cdot(Q_i)^{k_i-\delta_i}\cdot
\prod\limits_{\begin{array}{c} 0\leq j\leq s\\ j\not=i\end{array}}
\left(\begin{array}{c} k_j+\delta_j\\
\delta_j\end{array}\right)\cdot (Q_i-Q_j)^{-\delta_j} \\
 \le (n \cdot (2p^{2}+1))^{\sum_{0}^{s}\delta_{j}}
  \le  (2p^{2}+1)^{p} \cdot n^{p}
\end{eqnarray*}
where we have used that $\sum_{j=0}^{s}\delta_{j} = k_{i} \le p$.
 The total number of terms in the sum (\ref{fi1}) is
\begin{eqnarray}\label{partitions}
|P(s,k_i)| = \left(\begin{array}{c} s+k_i\\
k_i\end{array}\right) \leq n^{p}
\end{eqnarray}
since $s+ k_{i} \le n$ and $k_i\leq p$. This proves
(\ref{estim}).

Now consider the fractions $\frac{Q_i}{Q_i-Q_j}$ which appear in
(\ref{vol-formula}). For $0 \le j\leq m$, $j\not=i$ we have by
(\ref{lowerbound})
\begin{eqnarray}
\label{factors3} \left|\frac{Q_i}{Q_i-Q_j}\right|\leq
\frac{1}{\left| Q_i-Q_j\right|}  \leq 2p^2+1.
\end{eqnarray}

Note that $Q_j< Q_{m+1}\leq 0$ for any $m+2\leq j\leq s$. Using
statement (c) of Lemma \ref{lm2} we have for $m+2\leq j \le s$
\begin{eqnarray}\label{factors}
 0 < \frac{Q_i}{Q_i-Q_j}
\leq \frac{Q_i}{Q_i+ (2p+1)^{-1}}\leq \frac{1}{1+
(2p+1)^{-1}}=\frac{2p+1}{2p+2}.\end{eqnarray} If $Q_{m+1}<0$ then
estimate (\ref{factors}) continues to hold whereas if $Q_{m+1}=0$
then the corresponding factor equals $1$ and by Lemma \ref{lm2},
(e), using $k_{m+1}+1\leq p$, we get
\begin{eqnarray}\label{factors1}
\prod\limits_{\begin{array}{c} m < j\leq s \end{array}} \left|
\frac {Q_i}{Q_i-Q_j} \right| ^{k_j+1} \le a_0^{(\sum_{j = m+1}^{s}
(k_{j}+1)-p)} \le C' \cdot a_0^{n}
\end{eqnarray}
where the constant $C'$ depends on $p$ only. The number of
summands in formula (\ref{vol-formula}) equals $m+1$; by
(\ref{mki}) it is bounded above by $2p$.

Combining inequalities (\ref{estim}), (\ref{factors3}) and
(\ref{factors1}) we obtain an estimate of the form $1-r_J<
Cn^\gamma a_0^n$ where the constants $C$ and $\gamma$ depend on
$p$ but are independent of $n$. This clearly gives (\ref{estim1}).
\end{proof}

The following statement is equivalent to Proposition \ref{estim2};
we will need it in the proof of Theorem \ref{thm1}.

\begin{Proposition}\label{estim3} Given an integer $p\geq 1$,
there exist constants $C>0$ and $0<a<1$ such that for any $n\geq
p$ and any subset $J\subset \{1, 2, \dots, n\}$ of cardinality
$n-p$ one has
\begin{eqnarray}\label{estim4}
r_J < C\cdot a^n.
\end{eqnarray}
\end{Proposition}
\begin{proof}
The claimed statement follows from Proposition \ref{estim2} by
observing that $r_J = 1-r_{\bar J}$ where $\bar J$ denotes the
complement of $J$ in $\{1, \dots, n\}$.
\end{proof}


\begin{proof}[Proof of Theorem \ref{thm1} for $\mu=\mu_a$]
By Propositions \ref{prop2} and \ref{prop3} we have
$$b_p(n,\mu) = \sum_J r_J$$
where $J\subset \{1, \dots, n\}$ runs over all subsets containing
$1$ and being of cardinality either $p+1$ or $n-2-p$. By
Proposition \ref{estim2} each $r_J$ with $|J|=p+1$ contributes to
$b_p(n,\mu)$ a quantity exponentially close to $1$ and by
Proposition \ref{estim3} each term $r_J$ with $|J|=n-2-p$ is
exponentially small. Adding up all these contributions we arrive
at the desired inequality (\ref{ineqthm1}). \end{proof}

\section{Proof of Theorem \ref{thm1} for $\mu=\mu_b$}

The proof of Theorem \ref{thm1} in the case $\mu=\mu_b$ is quite
similar. Propositions \ref{estim2} and \ref{estim3} remain true
but their proofs are slightly different. The difference between
the two cases stems only from different simplices involved: for
$\mu=\mu_a$ we consider the simplex $A$ with vertices $c_0, \dots,
c_n$ and for $\mu=\mu_b$ we have to consider instead the simplex
$B$ with vertices $c'_0, \dots, c'_n$, see (\ref{verticesB}).

Let us examine the arguments of the proof of Proposition
\ref{estim2} when $c_i$ is replaced by $c'_i$. Inequality
(\ref{mki}) follows from (b) of Lemma \ref{lmqprime} and
inequality (\ref{ineq3}) follows from statement (e) of Lemma
\ref{lmqprime}. One introduces the points $Q'_0>Q'_1>\dots >Q'_s$
as the distinct values appearing in the sequence $q'_0, q'_1,
\dots, q'_n$. Instead of (\ref{lowerbound}) we have a simpler
inequality $|Q'_i-Q'_j|\geq 1$ where $i\not=j$ which is a
consequence of Lemma \ref{lmqprime}, (c).

Let us assume that $Q'_i>0$ for $i=0, \dots, m$ and $Q'_i\leq 0$
for $i=m+1, \dots, s$. We claim that for any $0\leq i\leq m$ the
quantity $F_i$ given by (\ref{fi1}) satisfies inequality
(\ref{estim}). Indeed, $|Q'_i|\leq 2p$ for $0\leq i\leq m$ (see
(a) and (b) of Lemma \ref{lmqprime}) and hence
$$|F_i|\leq \sum_{\delta\in P(s, k_i)} (2p)^{k_i-\delta_i}\cdot
n^{\sum \delta_j} \leq (2p)^p\cdot n^{2p} \leq (2p^2+1)^pn^{2p}.
$$
Here we have used inequalities (\ref{binom1}) and
(\ref{partitions}). To estimate (\ref{vol-formula}) from above we
note that for $0\leq i\leq m$ and $j\not=i$ one has
$$
\left| \frac{Q'_i}{Q'_i-Q'_j}\right| \, \leq\,  \left\{
\begin{array}{ll}
|Q'_i|\leq 2p & \mbox{for}\, \, j\leq m,\\ \\
1 & \mbox{for}\, \, j= m+1,\\ \\
\frac{Q'_i}{Q'_i+1} \leq \frac{2p}{2p+1} & \mbox{for} \, \, j \geq
m+2.
\end{array} \right.
$$
Here we have used statement (d) of Lemma \ref{lmqprime}. Combining
the obtained inequalities we obtain that the statement analogous
to Proposition \ref{estim2} holds for $\mu=\mu_b$. The remaining
arguments of the proof of Theorem \ref{thm1} for $\mu=\mu_b$ are
very similar to those described in the case $\mu=\mu_a$.

\end{document}